\documentclass[12pt,draft,amsmath,amscd,amsbsy,amssymb]{amsart}
\newtheorem{thm}{Theorem}[section]

\newtheorem{lem}[thm]{Lemma}
\newtheorem{prop}[thm]{Proposition}

\theoremstyle{definition}

\newtheorem{que}[thm]{Question}
\theoremstyle{remark}

\numberwithin{equation}{section}

\newcommand{\R}{\mathbb R}

\newcommand{\Z}{\mathbb Z}

\newcommand{\N}{\mathbb N}

\DeclareMathSymbol{\emptyset}{\mathord}{AMSb}{"3F}
\DeclareMathSymbol{\preccurlyeq}{\mathrel}{AMSa}{"34}

\newcommand{\lpm}{{\mathcal P}{\mathcal M}_{{\mathrm L}{\mathrm i} {\mathrm p}}}

\newcommand{\lm}{{\mathcal M}_{{\mathrm L}{\mathrm i} {\mathrm p}}}

\newcommand{\lpu}{{\mathcal U}{\mathcal M}_{{\mathrm L}{\mathrm i} {\mathrm p}}}

\newcommand{\lpf}{{\mathcal C}_{{\mathrm L}{\mathrm i} {\mathrm p}}}
\usepackage[all]{xy}
\begin{document}

\title[]{Nonexistence of linear operators extending Lipschitz (pseudo)metrics}
\author{ Michael Zarichnyi}
\address{Department of Mechanics and Mathematics,
Lviv National University, Universytetska Str. 1, 79000 Lviv,
Ukraine}
\email{topology@franko.lviv.ua, mzar@litech.lviv.ua}

%\thanks{}
\date{}
\dedicatory{}
\commby{}

%%% ----------------------------------------------------------------------

\begin{abstract} We provide an example of a zero-dimensional compact metric
space $X$ and its closed subspace $A$ such that there is no
continuous linear extension operator for the Lipschitz
pseudometrics on $A$ to the Lipschitz pseudometrics on $X$. The
construction is based on results of A. Brudnyi and Yu. Brudnyi
concerning linear extension operators for Lipschitz functions.
\end{abstract}

%%% ----------------------------------------------------------------------
\maketitle
%%% ----------------------------------------------------------------------

\section{Introduction}

The problem of extensions of metrics has a long history. It was F.
Hausdorff who first proved that any continuous metric defined on a
closed subset of a metrizable space can be extended to a continuous
metric defined on the whole space. C. Bessaga first considered the
problem of existence of linear extension operators for metrics (see
\cite{be2} and \cite{be1}) and provided a partial solution of this
problem. The problem was completely solved by T. Banakh \cite{b}
(see also a short proof in \cite{z}).

The problem of extension of Lipschitz and uniform (pseudo)metrics
has been considered in \cite{L}. It is known that any Lipschitz
pseudometric defined on a closed subset of a metric space admits an
extension which is a Lipschitz pseudometric defined on the whole
space. In this note we consider a problem of existence of linear
extension operators for Lipschitz pseudometrics. Up to the author's
knowledge, no affirmative results are obtained in this direction.
Our aim here is to construct a counterexample: there exists a
subset of a zero-dimensional compact metric space for which there
is no such an extension operator. The example is based on the
results from \cite{BB} concerning the linear extension operators
for Lipschitz functions.

Note that conditions for existence of extensions of continuous
functions are often equivalent to those for existence of extensions
of continuous pseudometrics (see, e.g., \cite{ASh}, \cite{A},
\cite{Sh}, \cite{Se}). It turns out that, in the case of linear
extension operators, one is able to proceed at least in one
direction, from pseudometrics to functions.

\section{Preliminaries}

Let $(X,d)$ be a compact metric space. Given a subset $A$ of $X$,
we say that a pseudometric $\varrho$ on $A$ is said to be {\it
Lipschitz} if there is $C>0$ such that $d(x,y)\le C\varrho(x,y)$,
for any $x,y\in A$. Also, a function $f\colon A\to\R$ is {\it
Lipschitz} if if there is $C>0$ such that $|f(x)-f(y)|\le Cd(x,y)$,
for every $x,y\in A$. Denote by $\lpm(A)$ (respectively $\lm(A)$,
$\lpf(A)$) the set of all Lipschitz pseudometrics (respectively
metrics, functions) on $A$. The set $\lpm(A)$ (respectively
$\lpf(A)$) is a cone (respectively linear space) with respect to
the operations of pointwise addition and multiplication by scalar.
We endow $\lpm(A)$ with the norm $\|\cdot\|_{\lpm(A)}$,
$$\|d\|_{\lpm(A)}=\sup\left\{\frac{d(x,y)}{\varrho(x,y)}\mid x\neq y\right\}$$
and $\lpf(A)$ with the seminorm $\|\cdot\|_{\lpf(A)}$,
$$\|f\|_{\lpf(A)}=\sup\left\{\frac{|f(x)-f(y)|}{\varrho(x,y)}\mid x\neq y\right\}$$
(in the sequel, we abbreviate $\|\cdot\|_{\lpm(A)}$ and
$\|\cdot\|_{\lpf(A)}$ to $\|\cdot\|_A$).

We say that a map $u\colon \lpm(A)\to\lpm(X)$ is an {\em extension
operator} for Lipschitz pseudometrics if the following holds:
\begin{enumerate}
  \item $u$ is linear (i.e. $u(d_1+d_2)=u(d_1)+u(d_2)$, $u(\lambda d)=\lambda
  u(d)$ for every $d_1,d_2\in \lpm(A)$, $\lambda\in\R_+$);

  \item $u(d)|(A\times A)=d$, for every $d\in\lpm(A)$;

  \item $u$ is continuous in the sense that $\|u\|=\sup\{\|u(d)\|_X\mid \|d\|_A\le1\}$
  is finite.
\end{enumerate}

This definition is a natural counterpart of those introduced in
\cite{BB} for the extensions of Lipschitz functions. The following
notation is introduced in \cite{BB}: $$\lambda(A,X)=\inf\{\|u\|\mid
u\text{ is a linear extension operator from }\lpf(A)\text{ to
}\lpf(X)\}.$$

Similarly, we put
\begin{align*}
 \Lambda(A,X)= & \inf\{\|u\|\mid u\text{ is a
linear extension operator from } \\
   & \lpm(A)\text{ to }\lpm(X)\}.
\end{align*}

It can be easily proved (cf. \cite{BB}) that
$\Lambda(X)=\sup\{\Lambda(A,X)\mid A\subset X\}$ is a bi-Lipschitz
invariant of $X$.

\section{Auxiliary results}
Given a metric space $X=(X,d)$ and $c>0$, we denote by $cX$ the
space $(X,cd)$.

\begin{lem}
For any metric space $(X,d)$, subset $A\subset X$, and $c>0$, we
have $\lambda(cS,cX)=\lambda(S,X)$.
\end{lem}
\begin{proof}
Let $\varphi\in\lpf(cS)$ and $\|\varphi\|_{cS}=K$. Then $\varphi$
can be also considered as an element of $\lpf(S)$ with
$\|\varphi\|_{S}=K/c$. There exists an extension $\bar\varphi\colon
X\to\R$ of $\varphi$ with $\|\bar\varphi\|_{X}\le
(K\lambda(S,X))/c$. Considering $\bar\varphi$ as an element of
$\lpf(cX)$, we see that $\|\bar\varphi\|_{cX}\le (K\lambda(S,X))$.
Therefore $\lambda(cS,cX)\le\lambda(S,X)$. Arguing similarly, one
can prove te opposite inequality.
\end{proof}
\begin{lem}\label{l:2}
Let a metric pair $(S_1,X_1)$ be a retract of a metric pair $(S,X)$
under a 1-Lipschitz retraction. Then $\lambda(S_1,X_1)\le
\lambda(S,X)$.
\end{lem}
\begin{proof}
 Let $r\colon X\to X_1$ be a 1-Lipschitz retraction such that
 $r(S)=S_1$. Given a Lipschitz function $f\colon S_1\to\R$, we see
 that $f\circ(r|S)$ is a  Lipschitz function on $S$ with
 $\|f\circ(r|S)\|_S=\|f\|_{S_1}$. There is an extension $g\colon X\to\R$ of $f\circ(r|S)$
 with $\|g\|_X\le \lambda(S,X)\|f\circ(r|S)\|_S$. Then $g|X_1$  is an extension
 of $f$ over $X_1$ with $\|g|X_1\|_{X_1}\le \lambda(S,X)\|f\|_{S_1}$.
\end{proof}

\begin{prop}\label{p:1} Let $S$ be a closed subset of a compact metric space $X$ with
$|S|\ge2$. The  following are equivalent:
\begin{enumerate}
  \item  there  exists a continuous linear extension operator
from $\lpm(S)$ to $\lpm(X)$;

  \item  there  exists a continuous linear  extension operator
from $\lm(S)$ to \\ $\lm(X)$.
\end{enumerate}

\end{prop}

\begin{proof} (1)$\Rightarrow$(2). Let $u\colon \lpm(S)\to \lpm(X)$ be
a continuous linear extension operator.

Note first that there exists $\tilde\varrho\in\lpm(X)$ such that
$\tilde\varrho(x,y)=0$ if and only if $(x,y)\in(S\times S)\cup
\Delta_X$ (by $\Delta_X $ we denote the diagonal of $X$). In order
to construct $\tilde\varrho$, for any $x,y\in X\setminus S$ with
$x\neq y$, consider the pseudometric $\varrho_{xy}$ on
$S\cup\{x,y\}$ defined by $$\varrho_{xy}|(S\times S)=0,\
\varrho_{xy}(x,y)=\varrho_{xy}(x,s)=\varrho_{xy}(y,s)=0$$ for any $s\in
S$. Denoting by $d$ the original metric on $X$ we see that
\begin{align*}
   & \varrho_{xy}(x,y)=1=(1/d(x,y))d(x,y), \\
   & \varrho_{xy}(x,s)=1\le(1/d(x,S))d(x,s)\\
   & \varrho_{xy}(y,s)=1\le(1/d(y,S))d(y,s)
\end{align*} for any $s\in S$. Therefore, $\varrho_{xy}$ is a
Lipschitz pseudometric with the Lipschitz constant $\max\{1/d(x,y),
1/d(x,S), 1/d(y,S)\}$. By the result of Luukkainen \cite{L}, there
exists a Lipschitz pseudometric, $\tilde\varrho_{xy}$, on $X$ which
is an extension of $\varrho_{xy}$.

There exist (necessarily disjoint) neighborhoods, $U_{xy}$ and
$V_{xy}$, of $x$ and $y$ respectively such that
$\tilde\varrho_{xy}(x',y')\neq0$, for every $x'\in U_{xy}$ and
$y'\in V_{xy}$. The family $$\{U_{xy}\times V_{xy}\mid x,y\in
((X\setminus S)\times (X\setminus S))\setminus\Delta_X\}$$ forms an
open cover of $((X\setminus S)\times (X\setminus
S))\setminus\Delta_X$ and, by separability of the latter set, there
exists a sequence $(x_i,y_i)$ in $((X\setminus S)\times (X\setminus
S))\setminus\Delta_X$ such that $$\bigcup\{U_{x_iy_i}\times
V_{x_iy_i}\mid i\in\N\}=((X\setminus S)\times (X\setminus
S))\setminus\Delta_X.$$ Let
$\displaystyle{\tilde\varrho=\sum_{i=1}^\infty
\frac{\tilde\varrho_{x_iy_i}}{2^i\|\tilde\varrho_{x_iy_i}\|_X}}$.
Then, obviously, $\tilde\varrho\in\lpm(X)$ is as required.

Now define an operator $\tilde u\colon \lm(S)\to\lm(X)$ as follows.
Let $x_0,y_0\in S$, $s_0\neq y_0$. Let $\tilde
u(\delta)=u(\delta)+\delta(x_0,y_0)\tilde\varrho$, for any
$\delta\in \lm(S)$. We leave to the reader an easy verification
that $\tilde u$ is a continuous extension operator.

(2)$\Rightarrow$(1). We are going to show that $\lm(S)$ is dense in
$\lpm(S)$. Let $\varepsilon>0$. Given $\varrho\in\lpm(S)$, we see
that $\varrho_1=\varrho+\varepsilon d'$, where $d'$ is the original
metric on $S$, is an element of $\lm(S)$ with
$\|\varrho-\varrho_1\|_S\le\varepsilon$.

Let $u\colon \lm(S)\to \lm(X)$ be a continuous linear extension
operator. Since $\lm(S)$ is dense in $\lpm(S)$ and the space
$\lpm(X)$ is complete, there exists a unique continuous extension,
$\tilde u\colon
\lpm(S)\to \lpm(X)$, of $u$. Obviously, $\tilde u$ is a continuous
linear extension operator.

\end{proof}

\section{Main result}

\begin{thm}\label{t:main}
There exists a closed subspace $A$ of a zero-dimensional compact
metric space $X$ for which there is no extension operator for
Lipschitz pseudometrics.
\end{thm}
\begin{proof}
We recall some results from \cite{BB}. Let $\Z_1^n(l)$ stand for
$\Z^n\cap[-l,l]^n$ endowed with the $\ell_1$-metric.

It is proved \cite[Lemma 10.5]{BB} that, for any natural $n$, there
exists $c_1>0$, $l(n)>0$, and a subset $Y_n\subset
\Z_1^n(l(n))$ such that
\begin{equation}\label{f:neq}
 \lambda(Y_n,\Z_1^n(l(n)))\ge
c_1\sqrt{n}.
\end{equation}
 Let $$X=\left(\coprod_{n=1}^\infty
\frac{1}{nl(n)}\Z_1^n(l(n))\right)\big/\sim,$$ where $\sim$ is the equivalence relation
which identifies all the origins, be the bouquet of the spaces
$\Z^n_1$, $n\in\N$. We naturally identify every
$\frac{1}{nl(n)}\Z_1^n(l(n)))$, $n\in\N$, with its copy, $X_n$, in
$X$. The space $X$ is endowed with the maximal metric, $\varrho$,
inducing the original metric on $X_n$, for every $n\in\N$.
Obviously, $X$ is a compact metric space. One can easily see that
$X$ is zero-dimensional.

By $\lpf(X,0)$ we denote the set of functions from $\lpf(X)$ that
vanish at $0\in X$ and by $\lpf^+(X,0)$ we denote the set of
nonnegative functions from $\lpf(X,0)$.

Let $S_n=\frac{1}{nl(n)}Y_n\subset X_n$, $n\in\N$. Suppose that
there exists an extension operator $u\colon\lpm(S)\to\lpm(X)$,
where $S=\coprod_{n=1}^\infty S_n\subset X$. Define a map $v\colon
\lpf(S)\to\lpf(X)$ by the following manner. First, let
$f\in\lpf^+(S,0)$, then the function
$$d_f\colon S\times S\to\R^+,\ d_f(x,y)=|f(x)-f(y)|,\ x,y\in X,$$ is a
Lipschitz pseudometric on $X$ and we let $v(f)(x)=u(d_f)(x,0)$.
Given $g\in\lpf(S,0)$, represent $g$ as the difference,
$g=g_1-g_2$, where $g_1,g_2\in \lpf^+(S,0)$ (say, $g_1=(|g|+g)/2$,
$g_2=(|g|-g)/2$), and let $v(g)=v(g_1)-v(g_2)$. Note that $v$ is
well-defined: if $g=g_1-g_2=g_1'-g_2'$, where $g_1,g_2,g_1',g_2'\in
\lpf^+(S,0)$, then $g_1+g_2'=g_1'+g_2$, whence, by the linearity of
$u$, we have
\begin{align*}
 v(g_1)(x)+v(g_2')(x)=&u(d_{g_1})(x,0)+u(d_{g_2'})(x,0)= u(d_{g_1}+d_{g_2'})(x,0)\\
  =&g_1(x)+g_2'(x)=g_1'(x)+g_2(x)=u(d_{g_1'}+d_{g_2})(x,0)\\ =&
  u(d_{g_1'})(x,0)+u(d_{g_2})(x,0)=v(g_1')(x)+v(g_2)(x),
\end{align*}
i.e., $v(g_1)+v(g_2')=v(g_1')+v(g_2)$.

 If $h\in\lpf(S)$, then $h-h(0)\in \lpf(S,0)$
and we put $v(h)=v(h-h(0))+h(0)$. By direct verification we show
that $v\colon\lpf(S)\to\lpf(X)$ is a linear extension operator with
$\|v\|<\infty$. Therefore, $\lambda(S,X)<\infty$.

For every $n$, denote by $r_n\colon X\to X_n$ the retraction that
sends the complement to $X_n$ to $0\in X_n$. Then, evidently, $r_n$
is a 1-Lipschitz retraction, $r_n(S)=S_n$ and, by Lemma \ref{l:2},
$\lambda(S_n, X_n)\le
\lambda(S,X)$. This obviously contradicts to
 inequality (\ref{f:neq}).
\end{proof}

It follows from Proposition \ref{p:1} that for the subset $S$ of
the space $X$ from the proof of Theorem \ref{t:main} there is no
continuous linear operator extending Lipschitz metrics.

\section{Remarks and open problems}

It easily follows from the proof of Theorem \ref{t:main} that
$\lambda(S,X)<\infty$, whenever $\Lambda(S,X)<\infty$, for any
subset $S$ of a metric space $X$.

\begin{que}\label{q:1}
Compare $\Lambda(S,X)$ and $\lambda(S,X)$.
\end{que}

We conjecture that $\Lambda(S,X)=\lambda(S,X)$; if, however, this
is not the case, one can ask for a pseudometric counterpart of any
result on linear extensions of Lipschitz functions. As an example,
we formulate the following question inspired by results from
\cite{bsh1}.

\begin{que}
 Let $(X,d)$ be a metric space and $\omega\colon\mathbb R_+\to\mathbb R_+$
 be a concave non-decreasing
function with $\omega(0)=0$. The function $d_{\omega}=\omega\circ
d$ is a metric on $X$. Are the properties
$\Lambda(X,d_{\omega})<\infty$ and $\Lambda(X,d)<\infty$
equivalent?
\end{que}

It is proved in \cite{tz} that there exists a linear operator which
extends partial pseudometrics with variable domain. The following
question is, in some sense, a strenthening of Question \ref{q:1}.
Given a compact metric space $X$, we let
$$\lpm=\bigcup\{\lpm(A)\mid A\text{ is a nonempty closed subset of
} X\}.$$ One can endow $\lpm$ with the following metric, $D$:
$$D(\varrho_1,\varrho_2)=\inf\{\|\tilde\varrho_1-\tilde\varrho_2\|_X
\mid \tilde\varrho_i\text{ is a Lipschitz extension of }\varrho_i\}.$$

\begin{que} Suppose that, for a metric space $X$, $\Lambda(X)<\infty$.
 Is there a continuous linear extension operator for partial
Lipschitz pseudometrics on $X$, i.e., a map $u\colon
\lpm\to\lpm(X)$ which is continuous with respect to the metric $D$
and whose restriction onto every $\lpm(A)$, where $A$ is a subset
of $X$, is linear?
\end{que}

Note that the question which corresponds to the above one for the
case of partial Lipschitz fuctions is also open; see \cite{ks} for
the results on simultaneous extensions of partial continuous
functions.

\begin{que}
The space $\lpm(X)$ can be endowed also with the topologies of the
uniform and pointwise convergence. Are there linear continuous
extensions operators from $\lpm(A)$ to $\lpm(X)$, where $A$ is a
subset of $X$, that are continuous in these topologies?
\end{que}

By $\lpu(Y)$ we denote the set of all Lipschitz ultrametrics on a
subset $Y$ of a zero-dimensional metric space. In \cite{tz1} (see
also \cite{s}) it is proved that there exists a continuous
extension operator that extends ultrametrics defined on a closed
subspace of a zero-dimensional compact metric space and preserves
the operation $\max$.

\begin{que}
Given a subset $A$ of a zero-dimensional metric space $X$, is there
a continuous extension operator for Lipschitz ultrapseudometrics,
$u\colon
\lpu(A)\to\lpu(X)$, that preserves the operation $\max$? is
homogeneous?
\end{que}

Apparently, the arguments of the proof of Theorem \ref{t:main} can
be applied to another situations in which there is no continuous
linear extension operator for a given class of functions; see,
e.g., \cite{bs}. We are going to return to these questions
elsewhere.

\end{document}